\documentclass[12pt]{article}
\usepackage[usenames]{color}
\usepackage[colorlinks=true,
linkcolor=webgreen,
filecolor=webbrown,
citecolor=webgreen]{hyperref}

\definecolor{webgreen}{rgb}{0,.5,0}
\definecolor{webbrown}{rgb}{.6,0,0}

\usepackage{amsmath,amssymb,psfig,epsfig,latexsym,graphicx}

\newtheorem{theorem}{Theorem}
\newtheorem{cor}[theorem]{Corollary}
\newtheorem{lem}[theorem]{Lemma}

\newtheorem{conj}{Conjecture}

\newcommand{\beql}[1]{\begin{equation}\label{#1}}
\newcommand{\eeq}{\end{equation}}
\newcommand{\bsq}{{\vrule height .9ex width .8ex depth -.1ex }}

\newcommand{\QQ}{{\mathbb Q}}
\newcommand{\ZZ}{{\mathbb Z}}
\newcommand{\RR}{\mathbb R}
\newcommand{\FF}{{\mathbb F}}
\DeclareMathOperator{\wt}{wt}
\DeclareMathOperator{\Fix}{Fix}
\DeclareMathOperator{\SL}{SL}
\DeclareMathOperator{\lcm}{lcm}
\newcommand{\Pn}{{\cal{P}}_n}
\newcommand{\RT }{{\cal{R}}}

\newcommand{\Zn}{\ensuremath{\ZZ[[x]]}}
\newcommand{\Zna}{\ensuremath{\ZZ[[x]]^{\ast}}}

\makeatletter
\def\@sect#1#2#3#4#5#6[#7]#8{\ifnum #2>\c@secnumdepth
     \def\@svsec{}\else
     \refstepcounter{#1}\edef\@svsec{\csname the#1\endcsname.\hskip .75em }\fi
     \@tempskipa #5\relax
      \ifdim \@tempskipa>\z@
        \begingroup #6\relax
          \@hangfrom{\hskip #3\relax\@svsec}{\interlinepenalty \@M #8\par}%
        \endgroup
       \csname #1mark\endcsname{#7}\addcontentsline
         {toc}{#1}{\ifnum #2>\c@secnumdepth \else
                      \protect\numberline{\csname the#1\endcsname}\fi
                    #7}\else
        \def\@svsechd{#6\hskip #3\@svsec #8\csname #1mark\endcsname
                      {#7}\addcontentsline
                           {toc}{#1}{\ifnum #2>\c@secnumdepth \else
                             \protect\numberline{\csname the#1\endcsname}\fi
                       #7}}\fi
     \@xsect{#5}}
\def\@begintheorem#1#2{\it \trivlist \item[\hskip \labelsep{\bf #1\ #2.}]}
\makeatother

\begin{document}
\begin{center}

{\large {\bf On the Integrality of $n$th Roots of Generating Functions }} \\
\vspace*{+.2in}
\end{center}

\begin{center}
Nadia Heninger${}^{(1)}$  \\
Computer Science Department \\
Princeton University \\
Princeton, NJ 08540 \\
(Email: nadiah@cs.princeton.edu)
\smallskip

E. M. Rains \\
Mathematics Department \\
University of California Davis \\
Davis, CA 95616 \\
(Email: rains@math.ucdavis.edu) \\
\smallskip

N.~J.~A.~Sloane${}^{(2)}$ \\
Internet and Network Systems Research Center \\
AT\&T Shannon Labs \\
180 Park Avenue,
Florham Park, NJ 07932--0971, USA \\
(Email: njas@research.att.com)

\end{center}

\smallskip
\begin{center}
DEDICATED TO THE MEMORY OF JACK VAN LINT (1932--2004).
\end{center}

\smallskip

\begin{center}
August 26, 2005; revised April 8, 2006
\end{center}

\smallskip

\begin{center}
{\bf Abstract}
\end{center}
Motivated by the discovery that the eighth root of the theta
series of the $E_8$ lattice
and the $24$th root of the theta
series of the Leech lattice 
both have integer coefficients,
we investigate the question of when an arbitrary element 
$f \in \RT$ (where $\RT = 1+x\ZZ [[x]]$)
can be written
as $f = g^n$ for $g \in \RT$, $n \ge 2$.  
Let $\Pn := \{ g^n \mid g \in \RT \}$ and let $\mu_n :=
n \prod_{p | n} p $.
We show among other things
that (i) for $f \in \RT , f \in \Pn \Leftrightarrow 
f \pmod{\mu_n} \in \Pn $, and
(ii) if $f \in \Pn$,
there is a unique $ g \in \Pn$
with coefficients mod $\mu_n/n$ such that $ f \equiv g^n$ (mod $\mu_n$).
In particular, if $f \equiv 1 ~(\mbox{mod}~ \mu_n)$ then $f \in \Pn$.
The latter assertion implies that the theta series of any
extremal even unimodular lattice in $\RR^n$ 
(e.g. $E_8$ in $\RR^8$)
is in $\Pn$ if $n$ is of the form $2^i3^j5^k \,(i \ge 3)$.
There do not seem to be any exact analogues for codes,
although we show that the weight enumerator of
the $r$th order Reed-Muller code
of length $2^m$ is in ${{\cal{P}}_{2^r}}$
(and similarly that the theta series of the
Barnes-Wall lattice $BW_{2^m}$ is in ${{\cal{P}}_{2^m}}$).
We give a number of other results and conjectures,
and establish a conjecture of Paul~D. Hanna that there
is a unique element $f \in \Pn$ ($n \ge 2$)
with coefficients restricted to the set
$\{ 1, 2, \ldots, n\}$.

\vspace{0.8\baselineskip}
${}^{(1)}$ Supported by the AT\&T Labs Fellowship Program. 
${}^{(2)}$ To whom correspondence should be addressed.

\vspace{0.8\baselineskip}

\vspace{0.8\baselineskip}

\section{Introduction}\label{SecIntro}

In June 2005, Michael Somos \cite{Somos05} observed that
the $12$-th root of the theta series of 
Nebe's extremal $3$-modular even lattice
in $24$ dimensions
(\cite{Nebe95}, \cite{Nebe98}, \cite{Quebbemann95},
sequence A004046 in \cite{OEIS})
appeared to have integer coefficients.
This led us to consider analogous questions
for other lattices, and we discovered that 
the cube root of the theta series of the $6$-dimensional lattice $E_6$,
the eighth root of the theta series of the $8$-dimensional lattice $E_8$,
and
the $24$th root of the theta series of the $24$-dimensional Leech
lattice $\Lambda _{24}$
also appeared to have integer coefficients.
Although it seemed unlikely (and still seems unlikely!)
that these results were not already known,
they were new to us, and so we considered the following 
general question.

Let $\Zn$ denote the ring of formal power series in $x$
with integer coefficients, let $\Zna$ denote the 
subset of $\Zn$ with constant term $\pm 1$ 
(that is, the set of units in $\Zn$),
and let $\RT  \subseteq \Zna$ be the elements with constant term $1$.
If $\Pn$ denotes the set $\{ g^n \mid g \in \RT  \}$,
when is a given $f \in \RT $ an element of $\Pn$ with $n \ge 2$?

In Section \ref{SecThms} we give some general conditions which
ensure that a series belongs to $\Pn$.
In Section \ref{SecTheta} we study the theta series of lattices 
and establish some general theorems which explain all the
above observations.  We also state some conjectures which would provide
converses to these theorems.
Section \ref{SecCodes} deals with the weight enumerators of codes.
Surprisingly (in view of the usual parallels between self-dual codes
and unimodular lattices, cf. \cite{SPLAG}, \cite{Elkies}, \cite{NRS06}),
there do not seem to be any exact analogues of the theorems
for theta series.
We show that the weight enumerator of
the $r$th order Reed-Muller code
of length $2^m$ is in ${{\cal{P}}_{2^r}}$ for $r=0, 1, \ldots, m$,
and make an analogous conjecture for extended BCH codes.
Similarly, we show that the theta series of the Barnes-Wall lattice
in $\RR ^{2^m}$ is in ${{\cal{P}}_{2^m}}$.
In Section \ref{SecSq} we consider the special case of series that
are squares, and report on a search for possible squares
in the {\em On-Line Encyclopedia of Integer Sequences} \cite{OEIS}.
This search led us to Paul Hanna's sequences, which are the subject of the final section.

It is worth mentioning that $\Zn$ is known to be a unique
factorization domain \cite{Samuel61},
although we will make no explicit use of this since we are concerned
only with the multiplicative group of units in $\Zn$.

\paragraph{Notation:}
If the formal power series $f(x) \in \Pn$ we will say that $f(x)$, or 
its sequence of coefficients, is ``an $n$th power''.
For a prime $p$, $|~|_p$ denotes the $p$-adic
valuation ($|0|_p := 0$; if $0 \ne r \in \QQ \, , \, 
r = p^a \, \frac{b}{c}
\mbox{~with~} a, b, c \in \ZZ, c \neq 0, \mbox{~and~}
\gcd(p,b) = \gcd(p,c) = 1,
\mbox{~then~} |r|_p := a$).
We will use the facts that
$|r!|_p < r/(p-1)$ for $r>0$,
$|\binom{p^i}{j}|_p = |p^i|_p-|j|_p$
(cf. \cite{Gouvea93}).

\section{ Conditions for $f$ to be an $n$th power }\label{SecThms}

We first show that, for investigating whether $f \in \RT $ is an $n$th power,
it is enough to consider $f$ mod $\mu_n$,
where 
$$
\mu_n := n \prod_{p | n} p \, .
$$

\begin{theorem}\label{ThTaylor}
For $f \in \RT , f \in \Pn$
if and only if
$f \pmod{ \mu_n} \in \Pn $.
\end{theorem}

\paragraph{Proof.}
We will show that, for $k \ge 1$,
the coefficients in $f^{1/n}$ are integers if and only if 
the coefficients in $(f + \mu_n x^k)^{1/n}$ are integers.
Let $\phi (f) := f^{1/n}$. By Taylor's theorem,
\begin{eqnarray}
\phi(f + \mu_n x^k)  & ~=~ & \sum_{r=0}^{\infty}  \,
\frac{ (\mu_nx^k)^r}{r!} \, \phi^{(r)}(f)  \nonumber \\
& ~=~ & \sum_{r=0}^{\infty}  \,
\frac{ (\mu_nx^k)^r}{r!} \,  r! \binom{\frac{1}{n}}{r} 
f^{ 1/n - r }
\nonumber \\
& ~=~ & f^{1/n}  \sum_{r=0}^{\infty} \mu_n^r \binom{\frac{1}{n}}{r} 
\frac{ x^{kr}}{f^r} ~. \nonumber 
\end{eqnarray}
Let $c := \mu_n^r \binom{\frac{1}{n}}{r}$.
For a prime $p$ dividing $n$,
$|c|_p = r|\mu_n|_p - r|n|_p - |r!|_p \ge 0$,
 by definition of $\mu_n$. For a prime
$p$ not dividing $n$, $1/n$ is a $p$-adic unit and again
 $|c|_p \ge 0$.
Hence $c \in \ZZ$. 
Since $f \in \RT $, $f^{-r}$ has integer coefficients,
and so
$(f + \mu_n x^k)^{1/n} = f^{1/n}g$ for some $g \in \RT $.
Thus the coefficients in $(f + \mu_n x^k)^{1/n}$ are integers
if and only if the coefficients in $f^{1/n}$ are integers.~~~$\bsq$

Since $1 \in \Pn$, we have:

\begin{cor}\label{CorMod}
If $f \in \RT $ satisfies $ f \equiv 1 \pmod{ \mu_n}$,
then $f \in \Pn $.
\end{cor}

\begin{cor}\label{Cor3}
Suppose $f = 1 + f_1 x + f_2 x^2 + \cdots \in \RT $. 
If $A$ and $B$ are positive integers such
that $\mu_n | AB$ and $\mu_n | A^2$,
then
$f(Ax) \in \Pn$ if $B | f_1$.
\end{cor}

This is an immediate consequence of Corollary \ref{CorMod}.
Similar conditions involving further coefficients
of $f$ can be obtained in the same way.

For example, if $n=2$, $f^{1/2}(4x)$ has
integer coefficients for any $f \in \RT $,
and $f^{1/2}(2x)$ has integer coefficients if
$2 | f_1$. (See Section \ref{SecSq} for more
about the case $n=2$.)

Furthermore, $n$th roots are unique mod $\mu_n/n$:

\begin{theorem}\label{ThMod2}
Given $f \in \Pn $, there is a unique
$g \in \RT $ mod $\mu_n/n$ such that $g^n \equiv f \pmod{\mu_n}$.
\end{theorem}

\paragraph{Proof.}
Given $f \in \Pn $, suppose $g \in \RT $ is such that
$g^n \equiv f \pmod{\mu_n}$.
We will show that,
for any $k \ge 1$,
$(g + \frac{\mu_n}{n} x^k)^n \equiv g^n \equiv f \bmod{\mu_n}$.
In fact,
$$
(g + \frac{\mu_n}{n} x^k)^n = g^n + \sum_{r=1}^{n} \binom{n}{r} \Big(\frac{\mu_n}{n}\Big)^r
x^{rk} g^{n-r} \, .
$$
Then for $r \ge 1$, $c := \binom{n}{r} \big(\frac{\mu_n}{n}\big)^r$
is divisible by $\mu_n$, because for primes $q$
not dividing $n$, $|c|_q = |\mu|_q = 0$,
while if $p$ divides $n$ then
$|c|_p \ge |n|_p - |r|_p +r
\ge |n|_p - |r|_p +p^{|r|_p}
\ge |\mu_n|_p = |n|_p +1$.
So we may reduce the coefficients of $g$ mod $\mu_n/n$.

Conversely, suppose $g^n \equiv h^n$ (mod $\mu_n$)
but $g \not\equiv h$ (mod $\mu_n/n$). Let $g$ and $h$
first differ at the $x^k$ term:
\begin{eqnarray}
g & = & 1 + g_1 x + \cdots + g_{k-1} x^{k-1} + \alpha x^k + \cdots ~,
\nonumber \\
h & = & 1 + g_1 x + \cdots + g_{k-1} x^{k-1} + \beta x^k + \cdots ~,
\nonumber
\end{eqnarray}
with $\alpha \not\equiv \beta$ mod $\mu_n/n$. 
Equating coefficients of $x^k$ in $g^n \equiv h^n$ (mod $\mu_n$)
gives
$n \alpha \equiv n \beta$ (mod $\mu_n$),
which implies $\alpha \equiv \beta$ (mod $\mu_n/n$),
a contradiction. So $g$ is unique.~~~$\bsq$

In the other direction, associated with any
$g \in (\ZZ/\frac{\mu_n}{n}\ZZ)[[x]]$ with constant term $1$
is a unique $f \in (\ZZ/\mu_n\ZZ)[[x]] \cap \Pn$,
namely $f := g^n \bmod{\mu_n}$.
So the elements of ${{\cal{P}}_2}$, for example, are enumerated by
infinite binary strings beginning with $1$.

We also note the following useful lemma.

\begin{lem}\label{LemmaI}
For $r, s \ge 1$,
$$
{{\cal{P}}_r} \cap {{\cal{P}}_s} = {{\cal{P}}_{\lcm (r,s)}}.
$$
\end{lem}

\paragraph{Proof.}
Clearly 
${{\cal{P}}_{\lcm (r,s)}} \subset {{\cal{P}}_r}, \, {{\cal{P}}_s}$.
On the other hand, suppose $f \in 
{{\cal{P}}_r} \cap {{\cal{P}}_s}$. Let $a,b$ be integers such that
$ar+bs = \gcd(r,s)$, and define
$$
g := (f^{\frac{1}{r}})^b (f^{\frac{1}{s}})^a \, .
$$
Then $g \in \RT$ and 
$g^{\lcm(r,s)} = g^{rs/\gcd(r,s)} = f$.~~~$\bsq$

\section{ Theta series of lattices }\label{SecTheta}

The theta series of an integral lattice $\Lambda$ in $\RR^d$
(that is, a lattice in which all inner products are integers)
is 
$$
\Theta_{\Lambda}(x) := \sum_{u \in \Lambda} x^{u \cdot u}  \in \RT  \, .
$$
The theta series of extremal lattices in various genera are
especially interesting in view of their connections with
modular forms and Diophantine equations
(\cite{SPLAG}, \cite{Scharlau99}, \cite{Serre88}).

\begin{lem}\label{LemmaAPF}
If $f\in 1+m x\ZZ[[x]]$ for some integer $m$, then for any integer $n$,
$$
f^n\in 1+m n' x\ZZ[[x]] \,,
$$
where $n' = \prod_{p|m} p^{|n|_p} ~(\mbox{or~} 0 \mbox{~if~} n=0)$.
\end{lem}

\paragraph{Proof.}
It suffices to consider the case $m=p^k$, $k>0$ and $n$ prime.  
If $n\ne p$, the claim is trivial, while otherwise, if $f=1+mg$, then
$$
(f^p-1)/m = \sum_{i=1}^{p} m^{i-1} \binom{p}{i} g^i \,.
$$
Every term on the right is a multiple of $p$,
and thus the claim follows.~~~$\bsq$

\begin{theorem}\label{th:E1}
If $\Lambda$ is an extremal even unimodular lattice in $\RR^d$,
$d$ a multiple of $8$,
then $\Theta_\Lambda (x) \in {\cal{P}}_n$, where $n$ is obtained from 
$d$ by discarding any prime factors other than $2$, $3$ and $5$.
\end{theorem}

\paragraph{Proof.}
Suppose $d = 8t = 2^i 3^j 5^k 7^\ell \cdots$ (with $i \ge 3$),
and let $a = \lfloor d/24 \rfloor = \lfloor t/3 \rfloor$.
Then $n = 2^i 3^j 5^k$ and $\mu_n$ is a divisor of $30n$.  
It is known that $\Theta_\Lambda (x)$ can be written in the form
\beql{Eq1}
\Theta_\Lambda (x) = \sum_{i=0}^{a} c_i \psi^{t-3i}(x)  \Delta^i(x) \,  ,
\eeq
where
\beql{Eq2a}
\psi(x)   :=  \Theta_{E_8}(x)  = 1 + 240 \sum_{m=1}^{\infty} \sigma_3(m) x^{2m}\,  ,
\eeq
\beql{Eq2b}
\Delta(x)   :=  x^2 \prod_{m=1}^{\infty} (1-x^{2m})^{24}\,  ,
\eeq
$\sigma_3(m)$ is the sum of the cubes of the divisors of $m$, 
and the coefficients $c_0 := 1, c_1, \ldots, c_a$ are such that
\beql{Eq4}
\Theta_{\Lambda} (x) = 1 + O(x^{2a+2}) \, .
\eeq
We will show that
\beql{Eq5}
\Theta_\Lambda(x) \equiv 1 \pmod{30n}\,  \, ,
\eeq
which by Corollary \ref{CorMod} implies the desired result.

We apply Lemma \ref{LemmaAPF},
taking $f=\psi, m=240, n=t, n'=2^{i-3} 3^j 5^k$, obtaining
$\psi^t(x) \equiv 1 \pmod{30n}$.  
By equating $\eqref{Eq1}$ and $\eqref{Eq4}$, 
we obtain an upper triangular system of equations for the $c_i$ 
with diagonal entries equal to $1$;
this implies inductively that for $i \ge 1$, $c_i \equiv 0 \pmod{30n}$,
and $\eqref{Eq5}$ follows.~~~$\bsq$

The theta series mentioned in Theorem \ref{th:E1}
is a modular form of weight $w=d/2 \equiv 0$ mod~$4$ for the full
modular group $\SL _2(\ZZ)$. More generally, we have:

\begin{theorem}\label{th:EMF}
Let $f(x)$ be the extremal modular form
of even weight $w$ for $\SL _2(\ZZ)$ (cf. \cite{Mallows75}).
Then  $f(x) \in {\cal{P}}_n$,
where $n$ is obtained from $2w$ by discarding
all primes $p$ such that $p-1$ does not divide $w$.
\end{theorem}

\paragraph{Proof.}
To show that the extremal modular form of weight $w$ is in 
${\cal P}_n$, it suffices to construct {\em any} modular form of weight $w$ 
congruent to $1 \mod \mu_n$; this form may even have denominators, as long 
as they are prime to $\mu_n$.
Indeed, the difference between such a form and the extremal form 
will be a cusp form with all leading coefficients a multiple of $\mu_n$; 
it follows as in the proof of Theorem \ref{th:E1} that such a cusp form has {\em 
all} coefficients a multiple of $\mu_n$.

In particular, one may consider the Eisenstein series.  
Every nonconstant coefficient of $E_w$ for $w$ even is a
multiple of $(-2w)/B_w$, where $B_w$ is a Bernoulli number, 
so it suffices to show that $\mu_n$ divides the denominator of $B_w/(2w)$.
By a result of Carmichael \cite{Carmichael},
$m$ divides this denominator 
if and only if the exponent of $\ZZ_m^*$ divides $w$.  
In particular, $2^{k+2}$ divides the denominator
if and only if $2^k$ divides $w$, 
while for odd primes, $p^{k+1}$ divides 
the denominator if and only if $p^k(p-1)$ divides $w$.  
The stated rule for $n$ follows.~~~$\bsq$

For $2$- and $3$-modular
lattices, we take powers of $\Theta_{D_4}$ and $\Theta_{A_2}$
respectively to determine $\mu_n$.  Presumably these results could be
improved by using the respective Eisenstein series instead.

\begin{theorem}\label{th:E2}
If $\Lambda$ is an extremal $2$-modular lattice in $\RR^d$,
$d$ a multiple of $4$,
then $\Theta_\Lambda (x) \in {\cal{P}}_n$, where $n$ is obtained from 
$d$ by discarding any prime factors other than $2$ and $3$.
\end{theorem}

\begin{theorem}\label{th:E3}
If $\Lambda$ is an extremal $3$-modular lattice in $\RR^d$,
$d$ a multiple of $2$,
then $\Theta_\Lambda (x) \in {\cal{P}}_{n/2}$, where $n$ is obtained from 
$d$ by discarding any prime factors other than $2$ and $3$.
\end{theorem}

It is a consequence of Theorems 
\ref{th:E1}, \ref{th:E2} and \ref{th:E3} that
that the theta series of the following
lattices are in ${\cal{P}}_d$, where $d$ (the subscript)
is the dimension of the lattice:
$D_4$ [sequence A004011 in \cite{OEIS}],
$E_8$ [A004009],
$BW_{16}$ [A008409],
$\Lambda_{24}$ [A008408]
and
Quebbemann's $Q_{32}$ [A002272].
Also, the theta series of the
Coxeter-Todd lattice $K_{12}$ [A004010]
is in ${\cal{P}}_6$,
and the theta series of Nebe's $24$-dimensional
lattice [A004006]
is in ${\cal{P}}_{12}$, establishing
Somos's conjecture mentioned in Section \ref{SecIntro}.
In the next section we will show more generally that the 
theta series of the Barnes-Wall lattice $BW_{2^m}$
is in ${{\cal{P}}_{2^m}}$ for all $m \ge 1$.

The coefficients of the $n$th roots in these
examples in general will not be the coefficients of
any modular form (at least, not in the sense of being
associated to any Fuchsian group). 
$\Theta_{E_8}(e^{2 \pi i z})$, for example,
has a zero in the open upper half-plane,
and so its eighth root has an algebraic singularity
in the upper half plane, and the coefficients have
exponential growth.

The coefficients of the $n$th roots 
also do not appear to have
any particular combinatorial significance.
For example, the theta series of the $D_4$ lattice is
$$
1+24\,{x}^{2}+24\,{x}^{4}+96\,{x}^{6}+24\,{x}^{8}+144\,{x}^{10}+96\,{x}^{12}+ \cdots \, ,
$$
in which the coefficient of $x^{2m}$ is the number
of ways of writing $2m$ as a sum of four squares,
while its fourth root [A108092] is
\begin{eqnarray*}
& 1 & + ~ 6\,{x}^{2}-48\,{x}^{4}+672\,{x}^{6}-10686\,{x}^{8}+185472\,{x}^{10}-3398304\,{x}^{12} \\
&&{} ~ +64606080\,{x}^{ 14} -1261584768\,{x}^{16}+25141699590\,{x}^{18}-509112525600\,{x}^{20}\\
&&{} ~ +10443131883360\,{x}^{22} -216500232587520\,{x}^{24}+4528450460408448\,{x}^{26} \\
&&{} ~-95438941858567104\,{x}^{28}+ 2024550297637849728\,{x} ^{30} - \cdots \, .
\end{eqnarray*}
Do these coefficients have any other interpretation?

\paragraph{Further examples.}
The extremal {\em odd} unimodular lattices
have been completely classified (cf. \cite{Conway78}, \cite[Chap. 19]{SPLAG}),
and the $\Pn$ to which their theta series belong are as follows:
$\Theta_{\ZZ^d} ~ (1 \le d \le 7) \in {\cal{P}}_d$,
$\Theta_{D_{12}^{+}}$ [A004533] $\in {\cal{P}}_4$,
$\Theta_{E_{7}^{2+}}$ [A004535] $\in {\cal{P}}_2$,
while the theta series of $A_{15}^{+}$ [A004536] and 
the odd Leech lattice  [A004537]
are only in ${\cal{P}}_1$.
This is a straightforward verification
since the theta series are known explicitly.

The theta series of both $E_6$ [A004007] and 
its dual $E_6^{\ast}$ [A005129]
are $\equiv 1 \bmod{9}$ (this follows from 
\cite[p. 127, Eqs. (121), (122)]{SPLAG}),
and so are in ${\cal{P}}_3$.

Michael Somos \cite{Somos05} has also pointed out that 
$x j(x) \in {\cal{P}}_{24}$, where
$j(x)$ is the modular function
$\frac{1}{x} + 744 + 196884 x + \cdots$ (\cite{Schoeneberg}).
This follows from $x j(x) = \psi(x)^3/\Delta(x)$.

We believe that the values of $n$ in Theorems \ref{th:E1} -- \ref{th:E3}
are best possible as far as the primes 2, 3, 5 and 7 are concerned.
For example, it is easy to check that
the theta series of the extremal even unimodular lattice in $\RR^{56}$
[A004673]
belongs to ${\cal{P}}_8$ but not ${\cal{P}}_{56}$.

The following conjecture also seems very plausible,
although again we do not have a proof:

\begin{conj}
Let $\Theta_{\Lambda}(x)$ be the theta series
of a $d$-dimensional lattice.
If $\Theta_{\Lambda}(x) \in \Pn$ then $n \le d$. 
$($In fact, we have not found any counterexample to
the stronger conjecture that $\Theta_{\Lambda}(x) \in \Pn$
implies that $n$ divides $d$.$)$
\end{conj}
Note that, considered as a formal power series,
$\Theta_{\Lambda}(x)$ determines the dimension $d$
(see \cite[p. 47, Eq. (42)]{SPLAG})---in Conway's terminology \cite{Conway97}, 
the dimension is an ``audible'' property.

\section{ Weight enumerators of codes }\label{SecCodes}
The weight enumerator of an $[n, \, k, \, d \, ]_q$ code (that is,
a linear code of length $n$, dimension $k$
and minimal Hamming distance $d$ over the field $\FF _q$) is
$$
W_C(x) := \sum_{c \in C} x^{\wt(c)} \, ,
$$
where $\wt$ denotes Hamming weight
(\cite{LintI}, \cite{MS77}).
Although the weight enumerators are polynomials,
the roots, if they exist, are normally infinite series.
There does not seem to be an analogue of Theorem  \ref{th:E1}
for extremal doubly-even binary self-dual codes, since
the weight enumerator of the $[24,12,8]_2$ Golay code,
$$
1 + 759 x^8 + 2576 x^{12} + 759 x^{16} + x^{24} \, ,
$$ 
is not in  ${\cal{P}}_n$ for any $n > 1$.
However, the weight enumerator of the $[8,4,4]_2$ Hamming code,
$1+14x^4+x^8$, is in ${\cal{P}}_2$
since it is congruent to $1+2x^4+x^8 \bmod{4}$,
although it is not in ${\cal{P}}_n$ for any $n > 2$.

This Hamming code is also the Reed-Muller code $RM(1,4)$
(cf. \cite{LintI}, \cite{MS77}).
More generally, we have:

\begin{theorem}\label{ThRM}
Let $W_{r,m}(x)$ denote the weight enumerator of the
$r$th order Reed-Muller code $RM(r,m)$, for $0 \le r \le m$,
and let $W_{r,m}(x) := W_{m,m}(x) = (1+x)^{2^m}$ for $r>m$.
Then for $r \le m$,
\beql{EqRM1}
W_{r,m}(x) \equiv ( 1 ~+~ x^{2^{m-r}} ) ^ {2^r} \, \pmod{2^{r+1}} \, ,
\eeq
and so by Theorem \ref{ThTaylor} is in ${\cal{P}}_{2^r}$.
\end{theorem}

We will deduce Theorem \ref{ThRM} from the following result:

\begin{theorem}\label{ThRM2}
For $0 \le r \le m+1$,
\beql{EqRM3}
W_{r,m+1}(x)-W_{r,m}(x^2) \equiv 0 \pmod{2^{m+1}} \, .
\eeq
\end{theorem}

\paragraph{Proof.}
Reed-Muller codes may be built up recursively from
\beql{EqRMu}
RM(r,m+1) = \{ (u,u+v) \mid
u \in RM(r,m), \, v \in RM(r-1,m) \} \, ,
\eeq
for $1 \le r \le m$,
with $RM(0,m+1) = \{ 0^{2^{m+1}}, 1^{2^{m+1}} \}$,
$RM(m+1,m+1) = \{ 0,1 \}^{2^{m+1}}$
(\cite[Chap. 13, Theorem 2]{MS77}).
Let $G$ be the group $(\FF_2^+)^{m+1}$
in its natural action on $C:=RM(r,m+1)$
(consisting of the diagonal action of $(\FF_2^+)^m$ on
$RM(r,m)$ and $RM(r-1,m)$ together with the involution swapping the two halves).
If $O(x)$ is the generating function for $G$-orbits, indexed by the
weight of the elements of the orbit, then by Burnside's Lemma,
$$
|G| \, O(x) = \sum_{g\in G} W_{\Fix_g(C)}(x) \, ,
$$
where $W_{\Fix_g(C)}(x)$ is the weight enumerator
of the subcode fixed by $g$.
For nonzero $g$, $W_{\Fix_g(C)}=W_{r,m}(x^2)$, from \eqref{EqRMu}.
Therefore
$$
|G| \, O(x)  = W_{r,m+1}(x) + (|G|-1) W_{r,m}(x^2) \, .
$$
Since $|G|=2^{m+1}$, the result follows immediately.~~~$\bsq$

Theorem \ref{ThRM} now follows from Theorem \ref{ThRM2}
by induction on $m$.
Another consequence of Theorem \ref{ThRM2} is:

\begin{cor} 
For any dyadic rational number $\lambda$
$($i.e., any element of $\ZZ[1/2])$
satisfying $0\le \lambda \le 1$,
and any integer $r\ge 0$, the sequence
\beql{EqRM4}
f_{r,m}(\lambda ) = |\{u\in RM(r,m) \mid \wt(u) = \lambda 2^m \}|
\, , \quad m = r, r+1, r+2, \ldots \, ,
\eeq
converges $2$-adically as $m\to\infty$.  
\end{cor}

Special cases of this were already known, but
in view of the many investigations of
weight enumerators of Reed-Muller codes
(\cite{Kasami76},
\cite[\S6.2]{LintL},
\cite{MS77},
\cite{Sloane70},
\cite{Sugino71},
\cite{Sugita96},
etc.),
it is worth putting the general remark on record.
For example, in the special case $\lambda = \frac{1}{2^r}$ it follows from
\cite[Chap. 13, Theorem 9]{MS77} that the limit in \eqref{EqRM4} is
$2^r / \prod_{i=1}^r (1-2^i)$.
Other special cases may be deduced from the results
in \cite{Sloane70} (or \cite[Chap. 15, Theorem 8]{MS77})
and \cite{Kasami76}.

The Nordstrom-Robinson, Kerdock and Preparata codes
are closely related to Reed-Muller codes
(\cite{Hammons94}, \cite{LintK}, \cite{MS77}).
The weight enumerator of the Nordstrom-Robinson code
of length $16$ is in ${\cal{P}}_2$, and more generally so is that
of the Kerdock code of length $4^m$, $m \ge 2$
(this follows immediately from \cite[Fig. 15.7]{MS77}).
It appears, although we do not have a proof,
that the weight enumerator of the Preparata code
of length $4^m$ is in ${\cal{P}}_{2^{m-3}}$.

There is a conjectural analogue of Theorem \ref{ThRM} for BCH codes:

\begin{conj}\label{ConjBCH}
Let $C$ be obtained by adding an overall parity check
to the primitive BCH code of length $2^m-1$
and designed distance $2t-1$, so that $C$
has length $n=2^m$ and minimal distance $d \ge 2t$.
We conjecture that the weight enumerator of $C$
is in ${\cal{P}}_{2^m/d'}$,
where $d'$ is the smallest power of $2 \ge 2t$.
\end{conj}

We have verified this for $m \le 6$.

Here are three further examples.
The Hamming weight enumerator of the $[12,6,6]_3$
ternary Golay code [A105683] is in ${\cal{P}}_4$,
and that of the $[18,9,8]_4$ extremal self-dual code $S_{18}$ over $\FF_4$
(\cite{Cheng90}, \cite{MacWilliams78}, A014487)
is in ${\cal{P}}_{18}$.
A more unlikely example is the weight enumerator of the $[48,23,8]_2$
Rao-Reddy code (\cite{Rao73}, \cite{MS77}, [A031137]),
\begin{eqnarray*}
1 & + & 7530\,{x}^{8}+92160\,{x}^{10}+1080384\,{x}^{12}
\nonumber \\
& & {} ~
+7342080\,{x}^{14}+34408911\,{x}^{16}+111507456\,{x}^{18}
\nonumber \\
& & {} ~
+255566784\,{x}^{20}+417404928\,{x}^{22}+492663180\,{x}^{24}
\nonumber \\
& & {} ~
+417404928\,{x}^{26}+255566784\,{x}^{28}+111507456\,{x}^{30}
\nonumber \\
& & {} ~
+34408911\,{x}^{32}+7342080\,{x}^{34}+1080384\,{x}^{36}
\nonumber \\
& & {} ~
+92160\,{x}^{38}+7530\,{x}^{40}+{x}^{48} \, ,
\nonumber 
\end{eqnarray*}
which is a square since it is congruent to $(1+x^8+x^{16}+x^{24})^2 \pmod{4}$.
(The square root is given in A108179.)

Barnes-Wall lattices are also closely related to Reed-Muller codes
\cite{SPLAG}, \cite{NRS02}, \cite{NRS06}.
It will be convenient here to normalize these lattices so that
the $2^m$-dimensional Barnes-Wall lattice $BW_{2^m}$
has minimal norm $2^{m-1}$ (making $BW_{2^m}$ a $2^{m-1}$-modular
lattice, cf. \cite{Quebbemann95}, in which all norms
are multiples of $2^{\lceil \frac{m}{2} \rceil } $).
Thus the first few instances are 
$$
BW_{2} = \ZZ ^2, \,
BW_{4} = D_4, \,
BW_{8} =\sqrt{2} \, E_8, \,
BW_{16} =\sqrt{2}\, \Lambda _{16}, \ldots \,.
$$
In particular, we see that for $m=1, \ldots, 4$,
$BW_{2^m}$ is in ${{\cal{P}}_{2^m}}$.
In fact, we have:
\begin{theorem}
The theta series of $BW_{2^m}$ in $\RR ^{2^m}$
is congruent to $1 \pmod{2^{m+1}}$ for $m \ge 1$,
and is thus in ${{\cal{P}}_{2^m}}$.
More precisely, for $m \ge 2$, we have
\beql{EqBWTh}
\frac{ \Theta _{BW_{2^m}} (x) - 1}{2^{m+1}}
~\equiv~
(1-2^{m-1}) \, \frac{\Theta _{BW_{2^{m-1}}}(x^2) - 1}{2^{m}}
\pmod {2^{m}} \, .
\eeq
\end{theorem}

\paragraph{Proof.}
For $m=1$, $BW_{2} = \ZZ ^2$ and
$\Theta _{BW_{2}} (x) \equiv 1 \pmod{4}$.
The automorphism group $G$ of $BW_{2^m}$
contains as a normal subgroup the extraspecial 
group $2_{+}^{1+2m}$ (cf. \cite{NRS02}).
For $m \ge 2$ the extraspecial group
consists of four conjugacy classes of $G$, with representatives, sizes and 
fixed sublattices as shown in Table \ref{TT1} (here $n = 2^m$):
\renewcommand{\arraystretch}{2.1}
\begin{table}[htb]
\caption{ }
$$
\begin{array}{ccc} 
\hline
\mbox{~Representative~} & \mbox{~Size~} & \mbox{~Fixed sublattice~} \\
\hline
1 & 1 & BW_{2^m} \\
-1 & 1 & 0 \\
\begin{bmatrix} 0 & I_n \\ I_n & 0 \end{bmatrix}  & 2^{2m} + 2^m -2 & \sqrt{2} \, BW_{2^{m-1}} \\
\begin{bmatrix} 0 & I_n \\ -I_n & 0 \end{bmatrix}  & 2^{2m} - 2^m  & 0  \\
\hline
\end{array}
$$
\label{TT1}
\end{table}
\renewcommand{\arraystretch}{1.0}

\noindent Then Burnside's  Lemma gives us the congruence
$$
\Theta _{BW_{2^m}} (x) + (2^{2m}-2^m+1) 
+ (2^{2m}+2^m-2) \, \Theta _{BW_{2^{m-1}}} (x^2)
~ \equiv ~ 0 \pmod{2^{2m+1}} \, ,
$$
which implies \eqref{EqBWTh}.~~~$\bsq$

This in particular implies that, for any dyadic rational
$\lambda \ge 1$, the coefficient of $x^{\lambda 2^{m-1}}$
(that is, the number of lattice vectors of norm equal to
$\lambda$ times the minimal norm) in
$$
\frac{ \Theta _{BW_{2^m}} (x) - 1}{2^{m+1}}
$$
converges to a $2$-adic limit.
For the kissing number itself, i.e. for $\lambda = 1$,
the limit is 
$\prod_{i=1}^\infty (1+2^i)$.

We end this section with a question:
Is there a simple way to test if a code has a
weight enumerator which is an $m$-th power?

\section{ Squares }\label{SecSq}

We know from Theorem \ref{ThTaylor} that to test if a given
$f(x) \in \RT$ is a square,
it is enough to consider $f(x) \bmod 4$,
and from Theorem \ref{ThMod2} that
if $f(x)$ {\em is} a square then
there is a unique binary series $g(x)$ associated with it.
There is a simple necessary and
sufficient condition for
$f(x)$ to be a square.

\begin{theorem}\label{ThSquare}
Given $f(x) := 1 + \sum_{r \ge 1} f_r x^r \in \RT$,
let $\bar{f}(x) := 1 + \sum_{r \ge 1} \bar{f_r} x^r$
be obtained by reducing the coefficients of $f(x)$ mod $4$.
If $\bar{f}_{2t} -g_t$ and $\bar{f}_{2t+1}$
are even for all $t \ge 0$,
where $g_0:=1, \, g_1, \ldots \in \ZZ/2\ZZ$ are
defined recursively by
\begin{eqnarray}\label{EqRec1}
\frac{\bar{f}_{2t} -g_t}{2} & \equiv & g_{2t} +
\sum_{r=1}^{t-1} g_r g_{2t-r} \pmod{2} \, , \nonumber \\
\frac{\bar{f}_{2t+1}}{2} & \equiv & g_{2t+1} +
\sum_{r=1}^{t} g_r g_{2t+1-r} \pmod{2} \, ,
\end{eqnarray}
then $f(x) \in {\cal{P}}_2$ and
\beql{EqRec2}
f(x) \equiv \bar{f}(x) \equiv g^2(x) := (1+\sum_{r=1}^{\infty} g_r x^r)^2 \pmod{4} \, .
\eeq
Conversely, if for some $t$ either $\bar{f}_{2t} -g_t$ or $\bar{f}_{2t+1}$
fails to be even, then $f(x) \notin {\cal{P}}_2$.
\end{theorem}

There is a simple {\em necessary} condition for $f(x)$ to
be a square, which generalizes to $p$th powers
for any prime $p$.

\begin{theorem}\label{ThSquare2}
Let $p$ be a prime.
If $f(x) := 1 + \sum_{r \ge 1} f_r x^r \in {\cal{P}}_p$,
say $f(x) = g(x)^p$, then
\beql{EqPP0}
f_r \equiv 0 \pmod{p} \mbox{ unless } p \mbox{ divides } r \, ,
\eeq
\beql{EqPP1}
g(x) \equiv 1 + f_p x + f_{2p} x^{2} + f_{3p} x^{3} + \cdots \pmod{p} 
\eeq
and
\beql{EqPP2}
f(x) \equiv (1 + f_p x + f_{2p} x^{2} + f_{3p} x^{3} + \cdots)^p \pmod{p^2} \, . 
\eeq
\end{theorem}

\paragraph{Proof.}
This follows immediately from Theorem \ref{ThMod2} and the fact that
$$
(1 + g_1 x + g_2 x^{2} + g_3 x^{3} + \cdots)^p
\equiv
1 + g_1 x^p + g_2 x^{2p} + g_3 x^{3p} + \cdots \pmod{p} \, .
$$
$\bsq$

The {\em On-Line Encyclopedia of Integer Sequences}
\cite{OEIS} is a database containing over $100,000$
number sequences.
We tested the corresponding formal power series to see which
were -- or at least appeared to be -- in ${\cal{P}}_2$.
As a first step we used the symbolic language {\em Maple}
\cite{Maple}
to weed out any series which did not begin $1 + \cdots$
or which had an obviously non-integral square root.
This produced $3030$ possible members of ${\cal{P}}_2$.
To reduce this number we discarded those series
which appeared to be congruent to 1 mod $4$, which left
$905$ candidates.

More detailed examination of these $905$ showed that most
of them could be grouped into one of the following 
(not necessarily disjoint) classes.

(1) Sequences which are obviously squares, usually with
a square generating function. 
These are often described as ``self-convolutions'' of other sequences.
For example, A008441, which gives the number of
ways of writing $n$ as the sum
of two triangular numbers,
with generating function $x^{-1/4}\eta(x^2)^4/\eta(x)^2$,
where $\eta(x)$ is the Dedekind eta function.

(2) Sequences which reduce mod $4$ to a square.
For example, periodic sequences of the form
$$
1,2,3, \ldots, k,
1,2,3, \ldots, k,
1,2,3, \ldots, k, \ldots \, ,
$$
are squares if and only if $k$ is a multiple of $4$.
More generally, any sequence which reduces mod $4$ to 
$1,2,3,4,5,6,\ldots $ is a square.

(3) Theta series of lattices and weight enumerators of codes,
as discussed in the preceding sections.

(4) McKay-Thompson series associated with 
conjugacy classes in the Monster simple group
(\cite{Ford94}, \cite{McKay90}, e.g. A101558).
As with the modular function $j(x)$ mentioned above,
the fact that these series are squares follows at once
from known properties.

(5) Sequences with an exponential generating function
involving trigonometric, inverse trigonometric, exponential, etc.,
functions.
One example out of many will serve as an illustration.
Vladeta Jovovi\'c's sequence A088313 \cite{Jovovic03}:
$$
1,2,7,36,241,1950,18271,193256,2270017,\ldots 
$$
gives the number of ``sets of lists'' with an odd number of
lists, that is, the number of partitions of $\{1, \ldots, n\}$ into an
odd number of ordered subsets (cf. Motzkin \cite{Motzkin71}).
There is no apparent reason why this should be a square. 
The analogous sequences for an even number of lists
(A088312) or with any number of lists (A000292) are not squares.
Jovovi\'c's sequence has exponential generating function
$$
\mbox{sinh}~\big(\frac{x}{1-x}\big) ~=~
x+\frac{2}{2!}{x}^{2}+{\frac {7}{3!}}{x}^{3}+
{\frac {36}{4!}}{x}^{4}+
{\frac {241}{5!}}{x}^{5}+{\frac{1950}{6!}}{x}^{6}+
{\frac { 18271}{7!}}{x}^{7}+ \, \cdots \, ,
$$
and is a square, since an elementary calculation shows that if
$$
\mbox{sinh}~\big(\frac{x}{1-x}\big) ~=~
\sum_{k=1}^{\infty} \, c_k \, \frac{x^k}{k!} \, ,
$$
then $c_k \equiv k \pmod{4}$.

(6) Paul Hanna's sequences, discussed in the following section.
These were the most interesting examples
that were turned up by our search. We were disappointed
not to find other sequences as challenging as these.

(7) Sequences whose square root proved to have a non-integral
coefficient once further terms were computed.

\section{ Paul Hanna's sequences }\label{SecHanna}

In May 2003, Paul D. Hanna \cite{Hanna2} contributed a family
of sequences to \cite{OEIS}.
For $k \ge 1$, the
$k$th Hanna sequence 
${H}_k := (1, h_1, h_2, \ldots)$ is
defined as follows: for all $n\ge 1$, $h_n$ is the smallest
number from the set $\{1, \ldots, k\}$ such that
$(1+h_1 x + h_2 x^2 + \cdots )^{1/k}$ has integer
coefficients. 
He asked if the sequences are well-defined and unique
for all $k$, and if they are eventually periodic.

For example, ${H}_2$ [A083952] is
$$
1,2,1,2,2,2,1,2,2,2,1,2,1,2,2,2,2,2,2,2,2,2,2,2,2,2,1,2,2,2,1,2,1,\ldots \, ,
$$
and the coefficients of its square root [A084202] are
$$
1,1,0,1,0,1,-1,2,-2,4,-6,10,-16,27,-44,75,-127,218,-375,650,-1130,
\ldots \, .
$$
The sequence ${H}_3$ [A083953] is
$$
1,3,3,1,3,3,3,3,3,3,3,3,1,3,3,2,3,3,2,3,3,1,3,3,2,3,3,3,3,3,2,3,3,3,3,
\ldots \, ,
$$
and the coefficients of its cube root [A084203] are
$$
1,1,0,0,1,-1,2,-2,2,0,-4,12,-24,38,-46,33,29,-176,443,-827,1222,-1310,
\ldots \, .
$$

\begin{theorem}\label{ThHanna0}
For all $k \ge 1$, ${H}_k$ is well-defined and is unique.
\end{theorem}

\paragraph{Proof.}
Suppose $f(x) := 1 + h_1 x + h_2 x^2 + \cdots = g(x)^k$,
where $g(x) := 1 + g_1 x + g_2 x^2 + \cdots$.
Then for $n \ge 1$, $h_n = k g_n + \Phi(g_1, \ldots, g_{n-1})$,
for some function $\Phi(g_1, \ldots, g_{n-1})$.
Write $\Phi(g_1, \ldots, g_{n-1}) = qk + r, ~ 0 \le r < k$.
If $r=0, h_n = k $ and $g_n = -(q-1)$,
while if $r>0$, $h_n=r$ and $g_n =-q$.~~~$\bsq$

We will analyze ${H}_2$  and ${H}_3$ in detail,
find generating functions for them,
and show that they are not periodic. 
We know from Section \ref{SecThms} that to study the $k$th root
$({H}_k)^{1/k}$ it is enough to look at its values mod $\mu_k/k$. 
The square root of  ${H}_2$ read mod $2$ gives
the binary sequence
$$
{S}_2 := (1,1,0,1,0,1,1,0,0,0,0,0,0,1,0,1,1,0,1,0,0,0,0,0,0,1,1,\ldots) 
$$
[A108336], and the cube root of ${H}_3$ read
mod $3$ gives 
$$
{S}_3 := (1,1,0,0,1,2,2,1,2,0,2,0,0,2,2,0,2,1,2,1,1,1,1,1,0,1,1,\ldots) 
$$
[A104405].

\begin{theorem}\label{ThHanna2}
The generating function $g(x) := 1 + x + x^3 + x^5 + x^6 + \cdots$ 
for ${S}_2$ satisfies $g(0)=1$ and
\beql{EqHan2}
g(x^2) + g(x)^2 \equiv \frac{2}{1-x} \pmod{4} \, .
\eeq
\end{theorem}

\paragraph{Proof.}
If $f(x)$ is the generating function for ${H}_2$,
we have $f(x) \equiv g(x)^2 \pmod{4}$. It follows (compare 
Theorem \ref{ThSquare}) that $f_{2t}=1$ if $g_t=1$, 
$f_{2t}=2$ if $g_t=0$, and $f_{2t+1}=2$.
Thus $f_{2t} \equiv 3 g_t + 2$ mod $4$. Hence
\beql{EqHGf2}
f(x) \equiv 3g(x^2) +\frac{2}{1-x^2} + \frac{2x}{1-x^2} \pmod{4} \, ,
\eeq
and \eqref{EqHan2} follows.~~~$\bsq$

\begin{cor}
${H}_2$ is not periodic.
\end{cor}

\paragraph{Proof.}
${H}_2$ is periodic if and only if ${S}_2$ is.
Suppose ${S}_2$ is periodic with period $\pi$.
Then $g(x) = p(x)/(1-x^{\pi})$,
where $p(x)$ is a polynomial of degree $\le \pi -1$.
 From \eqref{EqHan2},
\beql{EqHan2b}
\frac{p(x^2)}{1-x^{2 \pi}} +
\frac{p(x)^2}{(1-x^{\pi})^2} \equiv \frac{2}{1-x} \pmod{4} \, ,
\eeq
hence
$$
p(x^2)(1-x^{\pi}) + p(x)^2(1+x^{\pi}) \equiv
2 \frac{(1-x^{\pi})(1-x^{2\pi})}{1-x} \pmod{4} \, .
$$
The coefficient of $x^{3 \pi -1}$ is $0$ on the left,
$2$ on the right, a contradiction.~~~$\bsq$

Similar arguments apply to the ternary case; we
omit the details.

\begin{theorem}\label{ThHanna3}
Let $g(x) := 1 + x + x^4 + 2x^5 + 2x^6 + \cdots$ 
be the generating function for ${S}_3$,
and write it as $g(x) = g_{+}(x)+2g_{-}(x)$,
where $g_{+}(x)$ $($resp. $g_{-}(x))$ contains the
powers of $x$ with coefficient $1$ $($resp. $2)$. Then
$g(x)$ satisfies $g(0)=1$ and
\beql{EqHan3}
2g_{+}(x^3) + g_{-}(x^3)   + g(x)^3 \equiv \frac{3}{1-x} \pmod{9} \, .
\eeq
The generating function for ${H}_3$ is given by
\beql{EqHGf3}
f(x) \equiv \frac{3}{1-x} -2g_{+}(x^3) - g_{-}(x^3) \pmod{9} \, .
\eeq
\end{theorem}

\begin{cor}
${H}_3$ is not periodic.
\end{cor}

We have not studied the sequences ${H}_k$ for $k \ge 4$.

Another sequence of Hanna's is worth mentioning.
This is the sequence $a_0, a_1, a_2, \ldots$
defined by  $a_0=1$, and for $n>0$, $a_n$ is
the smallest positive number not already in the sequence
such that $(a_0 + a_1 x + a_2 x^2 + \cdots)^{1/3}$
has integer coefficients [A083349]:
$$
1,3,6,4,9,12,7,15,18,2,21,24,27,30,33,36,39,42,45,48,51,5,54,57,10,60,
\ldots \, .
$$
Although this sequence is similar in spirit to $H_3$,
there is no obvious relation between them.
Hanna  \cite{Hanna1} has shown that this sequence is a permutation 
of the positive integers.
No generating function is presently known.

\section*{ Postscript, Nov. 6, 2005 }
We cannot resist adding one further example, again a sequence
[A111983] studied by Paul Hanna. The series
$$
f(x) := \sum_{n=0}^{\infty} \, (2n+1) \, 8^n \, x^{ \frac{n(n+1)}{2} }
$$
is in  ${\cal{P}}_{12}$. Proof:
Mod 9, $f(x) \equiv  \sum_{0}^{\infty}
(-1)^n (2n+1) x^{n(n+1)/2} = \prod_{m=1}^{\infty} (1-x^m)^3$,
by an identity of Jacobi \cite[Th. 357]{Hardy},
so by Theorem \ref{ThTaylor} $f(x) \in {\cal{P}}_{3}$. 
Mod 8, $f(x) \equiv 1$, so $f(x) \in  {\cal{P}}_{4}$
by Corollary \ref{CorMod}, and then 
$f(x) \in  {\cal{P}}_{12}$ by Lemma \ref{LemmaI}.~~~$\bsq$

\section*{ Acknowledgments }

N.H. thanks the
AT\&T Labs Fellowship Program
for support during the summer of 2005 at Florham
Park, NJ, when this research was carried out. 
We thank Michael Somos for telling us about his discovery 
of the property of Nebe's lattice which prompted this work
and for further discussions about theta functions of lattices.
We also thank Andrew Granville for some helpful comments,
and Allan Wilks for some computations related to
Hanna's sequences.

\end{document}